(17/1/2018, 12h53)

# Note sur les lois locales conjointes de la fonction nombre de facteurs premiers

Gérald Tenenbaum

**Abstract.** Let $\alpha \in ]0,1]$ and let $Q_j$ $(1 \leqslant j \leqslant r)$ denote distinct irreducible polynomials with integer coefficients. We show that, for vectors with coordinates not exceeding a constant multiple of their mean, the joint local distribution of the number of prime factors of the $Q_j(n)$ for $x < n \leqslant x+x^\alpha$ is majorized by a constant multiple of the pairwise independency model, and we provide an upper bound for the constant in terms of the coefficients of the $Q_j$.

Notons $\omega(n)$ le nombre des facteurs premiers distincts d'un entier naturel $n$. Soit $\pi_{k,\ell}(x)$ le nombre des entiers $n$ n'excédant pas $x$ et tels que $\omega(n) = k$, $\omega(n+1) = \ell$. Un cas particulier d'un récent résultat de É. Goudout [2] stipule que la majoration

$$\pi_{k,\ell}(x) \ll \frac{x(\log_2 x)^{k+\ell-2}}{(\log x)^2(k-1)!(\ell-1)!} \qquad (x \geqslant 3) \tag{1}$$

est valable uniformément, pour toute constante $R > 0$, dans le domaine

$$1 \leqslant k, \ell \leqslant R \log_2 x.$$

Goudout établit en fait un théorème plus général, relatif à un produit de deux polynômes linéaires et fournit une borne uniforme dans les coefficients. Il mentionne également que sa démonstration s'adapte au cas d'un produit quelconque de polynômes linéaires.

Nous proposons ici une preuve simple d'une généralisation de (1), basée sur la méthode d'Erdős dans [1].

Nous considérons une famille $\{Q_j\}_{1\leqslant j\leqslant r}$ de polynômes irréductibles de $\mathbb{Z}[X]$, deux à deux premiers entre eux et sans diviseur fixe et posons $Q := \prod_{1\leqslant j\leqslant r} Q_j$. Nous désignons par $\varrho_j(m)$, resp. $\varrho_0(m)$, $(m \geqslant 1)$ le nombre de racines de $Q_j$, resp. de $Q$, dans $\mathbb{Z}/m\mathbb{Z}$, notons $D_j$ le discriminant de $Q_j$, $D$ celui de $Q$, et posons

$$g_j := \deg Q_j \quad (1 \leqslant j \leqslant r), \quad g := \sum_{1\leqslant j\leqslant r} g_j = \deg Q,$$
$$Q(X) = \sum_{0\leqslant i\leqslant g} \beta_i X^i, \quad \beta := \beta_g, \quad \|Q\| := \max_{1\leqslant i\leqslant g} |\beta_i|,$$
$$\varphi_j(n) := n\prod_{p|n}\Big(1 - \frac{\varrho_j(p)}{p}\Big) \quad (n \geqslant 1,\, 0 \leqslant j \leqslant r),$$

Ici et dans la suite nous réservons la lettre $p$ pour désigner un nombre premier.

En vertu du théorème des idéaux premiers, il existe des constantes positives $M_j$ $(1 \leqslant j \leqslant r)$ telles que

$$\left|\sum_{p\leqslant x} \frac{\varrho_j(p)}{p} - \log_2 x\right| \leqslant M_j \qquad (x \geqslant 2).$$

Notre résultat principal dépend fondamentalement de ces quantités.

Posons $M := \sum_{1\leqslant j\leqslant r} M_j$.

**Théorème 1.** *Soient $R > 0$, $r \geqslant 1$, $\alpha \in ]0,1[$. Nous avons uniformément pour $x \geqslant c_0\|Q\|$, $x^\alpha \leqslant y \leqslant x$, $k_j \in [1, R\log_2 x]$ $(1 \leqslant j \leqslant r)$,*

$$\sum_{\substack{x<n\leqslant x+y\\ \omega(Q_j(n))=k_j\ (1\leqslant j\leqslant r)}} 1 \ll \Big(\frac{\beta D}{\varphi_0(\beta D)}\Big)^K \frac{\mathrm{e}^M y}{(\log x)^r} \prod_{1\leqslant j\leqslant r} \frac{(\log_2 x + M_j)^{k_j-1}}{(k_j-1)!}, \tag{2}$$

*où $K$, $c_0$ et la constante implicite dépendent au plus de $\alpha$, $R$, $r$, $g$.*

*Remarques.* (i) Nous avons $\varrho_0(p) \leqslant g$ pour tout nombre premier $p$. Notant $\varphi$ la fonction indicatrice d'Euler, cela implique,

$$\frac{\beta D}{\varphi_0(\beta D)} \ll \Big(\frac{\beta D}{\varphi(\beta D)}\Big)^g \ll \Big(\log\{2+\omega(\beta D)\}\Big)^g.$$

(ii) Le membre de droite de (2) n'excède pas

$$\Big(\frac{\beta D}{\varphi_0(\beta D)}\Big)^K \frac{\mathrm{e}^{(R+1)M} y}{(\log x)^r} \prod_{1\leqslant j\leqslant r} \frac{(\log_2 x)^{k_j-1}}{(k_j-1)!}.$$

*Démonstration.* Soit $\varepsilon \in ]\alpha/4g, \alpha/3g[$. Pour chaque entier $n$ de $]x, x+y]$, désignons par $\xi_n$ le plus grand entier tel que

$$\prod_{\substack{p^\nu \| Q(n) \\ p \leqslant \xi_n}} p^\nu \leqslant x^{2g\varepsilon}.$$

Notant classiquement $P^+(n)$ (resp. $P^-(n)$) le plus grand (resp. le plus petit) facteur premier d'un entier $n$ avec la convention $P^+(1) = 1$, $P^-(1) = \infty$, nous obtenons la décomposition canonique

$$Q(n) = \prod_{1 \leqslant j \leqslant r} a_{jn} b_n$$

avec

$$P^+\Big(\prod_{1 \leqslant j \leqslant r} a_{jn}\Big) \leqslant \xi_n, \quad a_{jn} | Q_j(n) \ (1 \leqslant j \leqslant p_n r),$$
$$p_n := P^-(b_n) > \xi_n, \quad p_n^{v_n} \| b_n, \quad x^\varepsilon / p_n^{v_n} < \prod_{1 \leqslant j \leqslant r} a_{jn} \leqslant x^\varepsilon.$$

Désignons alors par $N_j(x)$ $(1 \leqslant j \leqslant 3)$ le nombre des entiers $n$ comptés dans le membre de gauche de (2) et satisfaisant respectivement aux conditions suivantes

$(N_1)$ $\qquad a_{1n} \cdots a_{rn} \leqslant x^{g\varepsilon}$ et $p_n > x^{\varepsilon/3}$,

$(N_2)$ $\qquad a_{1n} \cdots a_{rn} \leqslant x^{g\varepsilon}$ et $p_n \leqslant x^{\varepsilon/3}$,

$(N_3)$ $\qquad a_{1n} \cdots a_{rn} > x^{g\varepsilon}$.

Si $n$ est compté dans $N_1(x)$, alors $\omega(b_n) \leqslant E := 3(g+1)/\varepsilon$ dès que $x$ et $c_0$ sont assez grands. Notons $\mathcal{D}(U)$ le discriminant d'un polynôme $U \in \mathbb{Z}[X]$ et $\mathcal{R}(U, V)$ de résultant de $(U, V) \in \mathbb{Z}[X]^2$. L'argument développé dans [4], p. 134, une fois rectifié, comme indiqué dans [3], en tenant compte de la formule

$$c^{\deg U} d^{\deg V} \mathcal{D}(UV) = \mathcal{D}(U)\mathcal{D}(V)\mathcal{R}(U, V)^2,$$

valable pour tous polynômes $U$, $V$ de $\mathbb{Z}[X]$ de coefficients dominants respectifs $c$, $d$, nous permet alors d'écrire

$$a_{jn} = t_{jn} d_{jn} \quad (1 \leqslant j \leqslant r)$$

avec

$$t_{jn} | (\beta D)^\infty, \quad (d_{jn}, \beta D) = 1 \ (1 \leqslant j \leqslant r),$$
$$(d_{in}, d_{jn}) = 1 \ (1 \leqslant i < j \leqslant r).$$

Introduisons les noyaux sans facteur carré $t_{jn}^*$, $d_{jn}^*$ des $t_{jn}$, $d_{jn}$. Chaque entier $n$ compté dans $N_1(x)$ est donc tel que

$$\begin{cases} t_{jn}^* d_{jn}^* | Q_j(n) \quad (1 \leqslant j \leqslant r), \\ t_{1n}^* \cdots t_{rn}^* | (\beta D)^r, \\ \prod_{1\leqslant j\leqslant r} \mu\big(t_{jn}^*\big)^2 \mu\big(\prod_{1\leqslant j\leqslant r} d_{jn}^*\big)^2 = 1. \end{cases}$$

De plus, la condition $\omega\big(Q_j(n)\big) = k_j$ implique $\omega(t_{jn}d_{jn}) \in [k_j - E, k_j[$, où la borne supérieure découle de l'hypothèse $a_{jn} \leqslant x^{g\varepsilon}$.

En sommant selon les fibres

$$(t_{1n}^*, \ldots, t_{rn}^*, d_{1n}^*, \ldots, d_{rn}^*) = (t_1, \ldots, t_r, d_1, \ldots, d_r),$$

nous pouvons donc écrire

$$\begin{aligned} N_1(x) \leqslant & \sum_{\substack{(k_j - E)^+ \leqslant \kappa_j < k_j \\ (1\leqslant j\leqslant r)}} \sum_{t_1\cdots t_r | (\beta D)^r} \\ & \sum_{\substack{\omega(t_j d_j) = \kappa_j \ (1\leqslant j\leqslant r) \\ t_1 d_1 \cdots t_r d_r \leqslant x^{g\varepsilon} \\ (d_j, \beta D) = 1 \ (1\leqslant j\leqslant r)}} \prod_{1\leqslant j\leqslant r} \mu\big(t_j\big)^2 \mu\bigg(\prod_{1\leqslant j\leqslant r} d_j\bigg)^2 \sideset{}{^*}\sum_{x<n\leqslant x+y} 1 \end{aligned} \tag{3}$$

avec les conditions de sommation

$$(*) \begin{cases} t_1 \cdots t_r | (\beta^r D^r, Q(n)), \\ d_j | Q_j(n) \ (1 \leqslant j \leqslant r) \\ p | Q(n) \Rightarrow p | t_1 d_1 \cdots t_r d_r \text{ ou } p > x^{\varepsilon/3}. \end{cases}$$

Notant $T := t_1 \cdots t_r$, nous pouvons supposer sans perte de généralité que $\varrho_0(T) \prod_{1\leqslant j\leqslant r} \varrho_j(d_j) \geqslant 1$, puisque $N_1(x) = 0$ dans le cas contraire. La somme intérieure peut alors être majorée par le crible de Selberg comme au lemme 3.4 de [6] puisque l'ensemble

$$\mathcal{A} := \big\{n \in ]x, x+y] : T | Q(n),\, d_j | Q_j(n)\, (1 \leqslant j \leqslant r)\big\}$$

vérifie $|\mathcal{A}| = X + O(R)$ avec

$$X := y \frac{\varrho_0(T)}{T} \prod_{1\leqslant j\leqslant r} \frac{\varrho_j(d_j)}{d_j} \geqslant x^{2\alpha/3}, \quad R := \varrho_0(T) \prod_{1\leqslant j\leqslant r} \varrho_j(d_j) \leqslant x^{g\varepsilon} \leqslant x^{\alpha/3}.$$

Nous obtenons

$$\begin{aligned}
\sum_{x<n\leqslant x+y}{}^{*} 1 &\ll X \prod_{\substack{p\leqslant x^{\varepsilon/3}\\ p\nmid t_1d_1\cdots t_rd_r}} \Big(1-\frac{\varrho_0(p)}{p}\Big)\\
&\ll X \prod_{\substack{p\leqslant x^{\varepsilon/3}\\ p\nmid \beta Dd_1\cdots d_r}} \Big(1-\sum_{1\leqslant j\leqslant r}\frac{\varrho_j(p)}{p}\Big)\\
&\ll X\Big(\frac{\beta D}{\varphi_0(\beta D)}\Big)^r \prod_{1\leqslant j\leqslant r}\frac{d_j}{\varphi_j(d_j)}\exp\Big\{-\sum_{1\leqslant j\leqslant r}\sum_{p\leqslant x^{\varepsilon/3}}\frac{\varrho_j(p)}{p}\Big\}\\
&\ll \mathrm{e}^M\Big(\frac{\beta D}{\varphi_0(\beta D)}\Big)^r \prod_{1\leqslant j\leqslant r}\frac{\varrho_j(d_j)}{\varphi_j(d_j)}\frac{y\varrho_0(T)}{T(\log x)^r}\cdot
\end{aligned}$$

Notant $\tau_r(n)$ le nombre de décompositions d'un entier $n$ en produit de $r$ facteurs, il suit

$$N_1(x) \ll \frac{Hy}{(\log x)^r}\sum_{T\mid(\beta D)^r}\frac{\varrho_0(T)\tau_r(T)}{T}S(T),$$

avec

$$\begin{aligned}
H &:= \mathrm{e}^M\Big(\frac{\beta D}{\varphi_0(\beta D)}\Big)^r,\\
S(T) &:= \sum_{\substack{(k_j-E-\omega(T))^+\leqslant\kappa_j<k_j\\ (1\leqslant j\leqslant r)}}\prod_{1\leqslant j\leqslant r}\frac{1}{\kappa_j!}\Big(\sum_{p\leqslant x}\frac{\varrho_j(p)}{p-\varrho_j(p)}\Big)^{\kappa_j}\\
&\ll \big\{\omega(T)+E\big\}^r\prod_{1\leqslant j\leqslant r}\frac{(\log_2 x+M_j)^{k_j-1}}{(k_j-1)!}(R+1)^{r\omega(T)+rE},
\end{aligned}$$

compte tenu de la majoration

$$\begin{aligned}
\frac{1}{\kappa_j!}\Big(\sum_{p\leqslant x}\frac{\varrho_j(p)}{p-\varrho_j(p)}\Big)^{\kappa_j} &\ll \frac{1}{\kappa_j!}\Big(\log_2 x+M_j\Big)^{\kappa_j}\\
&\ll \frac{(\log_2 x+M_j)^{k_j-1}}{(k_j-1)!}\Big(1+\frac{k_j}{\log x}\Big)^{\ell}\leqslant (R+1)^{\ell}\frac{(\log_2 x+M_j)^{k_j-1}}{(k_j-1)!},
\end{aligned}$$

valable pour $k_j - \ell \leqslant \kappa_j < k_j \leqslant R\log_2 x$. Nous pouvons donc écrire

$$\begin{aligned} N_1(x) &\ll \frac{Hy}{(\log x)^r} \prod_{1\leqslant j\leqslant r} \frac{(\log_2 x + M_j)^{k_j-1}}{(k_j-1)!} \sum_{T|(\beta D)^r} \frac{\varrho_0(T)\tau_r(T)(R+2)^{r\omega(T)}}{T} \\ &\ll \Big(\frac{\beta D}{\varphi_0(\beta D)}\Big)^K \frac{\mathrm{e}^M y}{(\log x)^r} \prod_{1\leqslant j\leqslant r} \frac{(\log_2 x + M_j)^{k_j-1}}{(k_j-1)!}. \end{aligned}$$

Considérons à présent un entier $n$ compté dans $N_2(x)$. Nous avons $p_n^{v_n} > x^{g\varepsilon}$ et $p_n \leqslant x^{\varepsilon/3}$. Pour chaque nombre premier $p$ n'excedant pas $x^{\varepsilon/3}$, désignons par $\nu(p)$ le plus petit entier tel que $p^{\nu(p)} > x^{g\varepsilon}$. Alors $\nu(p) \geqslant 3g$ et $p^{\nu(p)-1} \leqslant x^{g\varepsilon}$, donc $p^{\nu(p)} \leqslant x^{2g\varepsilon} \leqslant y$. La majoration universelle $\varrho_0(p^\nu) \leqslant gp^{\nu(1-1/g)}$, établie par Stewart dans [5] pour toute puissance de nombre premier $p^\nu$, nous permet alors d'écrire

$$\begin{aligned} N_2(x) &\leqslant \sum_{p\leqslant x^{\varepsilon/3}} \sum_{\substack{x<n\leqslant x+y \\ Q(n)\equiv 0 \,(\mathrm{mod}\, p^{\nu(p)})}} 1 \leqslant 2 \sum_{x^{g\varepsilon/\nu(p)}<p\leqslant x^{\varepsilon/3}} \frac{y\varrho_0\big(p^{\nu(p)}\big)}{p^{\nu(p)}} \\ &\ll y \sum_{x^{g\varepsilon/\nu(p)}<p\leqslant x^{\varepsilon/3}} \frac{1}{p^{\nu(p)/g}} \ll \frac{y}{x^{2\varepsilon/3}}. \end{aligned}$$

Cette majoration est clairement compatible avec l'estimation annoncée (2).

Pour estimer $N_3(x)$, nous introduisons le paramètre $q_n := P^+(a_{1n}\cdots a_{rn})$ et notons que $\omega(b_n) \leqslant \eta(q_n) := (g+1)(\log x)/\log q_n$. Nous avons

$$\begin{aligned} N_3(x) \leqslant \sum_{q\leqslant x^{2g\varepsilon}} \sum_{\substack{(k_j-\eta(q))^+\leqslant\kappa_j<k_j \\ (1\leqslant j\leqslant r)}} \sum_{t_1\cdots t_r|(\beta D)^r} \\ \sum_{\substack{\omega(t_jd_j)=\kappa_j\ (1\leqslant j\leqslant r) \\ P^+(t_1d_1\cdots t_rd_r)=q \\ (d_i,\beta Dd_j)=1\ (1\leqslant i<j\leqslant r) \\ x^{g\varepsilon}<t_1d_1\cdots t_rd_r\leqslant x^{2g\varepsilon}}} \prod_{1\leqslant j\leqslant r} \mu(t_j)^2 \mu\Big(\prod_{1\leqslant j\leqslant r} d_j\Big)^2 \sum_{n\leqslant x}^{**} 1 \end{aligned}$$

avec les conditions

$$(**)\begin{cases} t_1\cdots t_r|(\beta^r D^r, Q(n)) \\ d_j|Q_j(n)\ (1\leqslant j\leqslant r) \\ p|Q(n) \Rightarrow p|t_1d_1\cdots t_rd_r \text{ ou } p>q. \end{cases}$$

En employant à nouveau le crible de Selberg, nous obtenons

$$\sum_{n\leqslant x}^{**} 1 \ll \frac{Hy}{(\log q)^r}\frac{\varrho_0(t_1\cdots t_r)}{t_1\cdots t_r}\prod_{1\leqslant j\leqslant r}\frac{\varrho_j(d_j)}{\varphi_j(d_j)}\cdot$$

Nous majorons alors la somme $r$-uple par la méthode de Rankin en introduisant un poids $(t_1d_1\cdots t_rd_r/x^{g\varepsilon})^v$, avec $v := C/\log q$, où $C$ est une constante assez grande. Il suit comme précédemment

$$N_3(x) \ll H\sum_{T\mid\beta^rD^r}\frac{\varrho_0(T)\tau_r(T)}{T^{1-v}}\sum_{q\leqslant x^{g\varepsilon}}\frac{y}{q(\log q)^r x^{Cg\varepsilon/\log q}}$$
$$\sum_{\substack{(k_j-\eta(q)-\omega(T))^+\leqslant\kappa_j<k_j\\(1\leqslant j\leqslant r)}}\prod_{1\leqslant j\leqslant r}\frac{(\log_2 q+M_j)^{\kappa_j}}{\kappa_j!}\cdot$$

Pour chaque $r$-uple $(\kappa_1,\ldots,\kappa_r)$ apparaissant dans la somme intérieure, nous avons

$$\prod_{1\leqslant j\leqslant r}\frac{(\log_2 q+M_j)^{\kappa_j}}{\kappa_j!}\leqslant\prod_{1\leqslant j\leqslant r}\frac{(\log_2 x+M_j)^{k_j-1}}{(k_j-1)!}(R+1)^{r\eta(q)+r\omega(T)}.$$

Dès que $Cg\varepsilon > 1+r(g+1)\log(R+1)$, nous obtenons donc, pour une constante $K$ convenable,

$$\begin{aligned}
&N_3(x)\\
&\ll\Big(\frac{\beta D}{\varphi_0(\beta D)}\Big)^K\frac{\mathrm{e}^M y}{(\log x)^r}\prod_{1\leqslant j\leqslant r}\frac{(\log_2 x+M_j)^{k_j-1}}{(k_j-1)!}\sum_{q\leqslant x^\varepsilon}\Big(\frac{\log x}{\log q}\Big)^{2r}\frac{x^{-1/\log q}}{q}\\
&\ll\Big(\frac{\beta D}{\varphi_0(\beta D)}\Big)^K\frac{\mathrm{e}^M y}{(\log x)^r}\prod_{1\leqslant j\leqslant r}\frac{(\log_2 x+M_j)^{k_j-1}}{(k_j-1)!}\cdot
\end{aligned}$$

□

**Remerciements.** L'auteur tient à exprimer ses chaleureux remerciements, d'une part, à l'arbitre pour avoir détecté une erreur significative dans une version initiale de ce travail et, d'autre part, à Régis de la Bretèche pour son aide amicale et pertinente.

## Bibliographie

Gérald Tenenbaum
Institut Élie Cartan
Université de Lorraine
BP 70239
54506 Vandœuvre Cedex
France

internet : gerald.tenenbaum@univ-lorraine.fr